\documentclass[12pt]{amsart}
\usepackage[dvips]{graphicx}
\usepackage{amssymb}
\usepackage{latexsym}

\def\C{\mathbb C}
\def\bC{\mathbf{\overline{C}}}
\def\N{\mathbb N}
\def\length{\operatorname{length}}
\def\area{\operatorname{area}}

\def\re{\operatorname{Re}}
\def\im{\operatorname{Im}}
\def\arg{\operatorname{arg}}
\newtheorem{la}{Lemma}
\newtheorem{theorem}{Theorem}
\newtheorem*{dca}{Denjoy-Carleman-Ahlfors Theorem}
\newtheorem*{ivth}{Iversen's Theorem}
\newtheorem*{cor}{Corollary}
\newtheorem{prop}{Proposition}
\theoremstyle{remark}
\newtheorem*{rem}{Remark}
\newtheorem*{ack}{Acknowledgment}
\begin{document}
\title[Direct singularities]{Direct singularities
and completely invariant domains
of entire functions}
\author{Walter 
Bergweiler}\thanks{The first author was supported by the Alexander von
Humboldt Foundation and by the G.I.F.,
the German--Israeli Foundation for Scientific Research and
Development, Grant G-809-234.6/2003} 
\address{Mathematisches Seminar,
Christian--Albrechts--Universit\"at zu Kiel,
Lude\-wig--Meyn--Str.~4,
D--24098 Kiel,
Germany}
\email{bergweiler@math.uni-kiel.de}
\author{Alexandre 
Eremenko}\thanks{The second author was 
supported by NSF grants DMS-0555279 and DMS-0244547}
\address{Department of Mathematics,
Purdue University, West Lafayette, IN 47907, USA}
\email{eremenko@math.purdue.edu}
\date{\today}
\begin{abstract}
Let $f$ be a transcendental entire function which omits
a point $a\in\C$. 
We show  that if $D$ is a simply
connected domain which does not contain $a$, then the full preimage
$f^{-1}(D)$ is disconnected. Thus, in dynamical context,
if an entire function has a completely invariant domain and omits
some value, then the omitted value belongs to the completely
invariant domain. We conjecture that the same
property holds if $a$ is a {\em locally omitted value}
(i.e., the projection of a direct singularity of $f^{-1}$).
We were able to prove this conjecture
for entire functions of finite order.
We include some auxilliary results on
singularities
of $f^{-1}$ for entire functions $f$, which can be of independent
interest.
\end{abstract}
\maketitle
\section{Introduction and results}
The question considered in this paper is motivated by
dynamics of entire functions~\cite{Be,EL}.
A component $D$ of the Fatou set of an entire function $f$ is called a
{\em completely invariant domain} if $f^{-1}(D)=D$. This is a stronger
property than simple invariance $f(D)\subset D$. 

In what follows, all entire functions are assumed to be transcendental.
It follows from a result of Baker~\cite[Theorem~1]{Bak75} that
all invariant components of the  Fatou set of such a function 
are simply connected. Baker also  proved
that at most one completely invariant domain
can exist~\cite{Ba}, and if $f$ has a completely
invariant domain, then all 
critical values (and thus all critical points) of $f$ are contained in
it~\cite[Theorem~2]{Bak75}.

In~\cite[Lemma~11]{EL}, the latter result of Baker was extended to the
logarithmic singularities of $f^{-1}$: a completely invariant domain
must contain all projections of
logarithmic singularities of~$f^{-1}$.
In this paper, we consider possibilities of extension
of this result to other types of singularities of $f^{-1}$.

A~point $a\in\C$ is {\em omitted} by an entire function $f$
if $f(z)\neq a$ for $z\in\C$.
A~point $a\in\C$ is {\em locally omitted} by $f$
if there exists $r>0$ and a component $G$ of the set $f^{-1}(B(a,r))$
such that $f(z)\neq a$ in~$G$. Here and in what follows, we
use the notation $B(a,r)$ for a disc of radius $r$ centered at $a\in\C$.
According to Iversen's classification of singularities, which 
will be recalled in Section~\ref{preliminaries}, a value is 
locally omitted if and only if it is the projection of a direct
singularity of~$f^{-1}$.
In the special case that $f:G\to B(a,r)\setminus\{a\}$ 
is a universal covering we say
that $a$ is the projection of a logarithmic singularity of~$f^{-1}$.

It is known that an omitted value does not have to be the projection of a
logarithmic singularity. An example of this is
\begin{equation}
\label{fexample}
f(z)=\exp\left(\sum_{k=1}^\infty
\left(\frac{z}{2^k}\right)^{2^k}\right).
\end{equation}
We will analyse this example in the end of the paper.

\begin{theorem}\label{t1}
Let $f$ be an entire transcendental function omitting
a point $a\in\C$, and let $D$ be a simply connected region
that does not contain~$a$. Then $f^{-1}(D)$ is disconnected.
\end{theorem}

\begin{cor}
Let $f$ be a transcendental entire function having a
completely invariant domain~$D$. If $f$ omits
a point then this point belongs to~$D$.
\end{cor}

We conjecture that Theorem~1 and the Corollary
can be extended to locally omitted values.
Paper~\cite{EL} contains a statement that the Corollary can be proved
for locally omitted values 
in the same way as for projections of logarithmic singularities.
However, the argument given in~\cite{EL} does not apply to
locally omitted values of arbitrary entire functions.
So the conjecture remains open.

In this paper we prove the conjecture for functions of finite order.
Namely we establish the following.

\begin{theorem}\label{t2} Let $f$ be an entire function of finite order,
and let $a\in\C$ be either a critical value or a locally omitted value.
If $D$ is a simply connected region that does not contain $a$,
then $f^{-1}(D)$ is disconnected.
\end{theorem}

Iversen's Theorem (stated in Section~2) implies that a
locally omitted value has to be an asymptotic value.
There is an example~\cite{Be2} of an entire function
of finite order with a completely invariant domain $D$
and an asymptotic value that does not belong to~$D$.

It is interesting that Theorem~2 has a converse:

\begin{theorem}\label{3}
Let $f$ be an entire function of finite order, and let $a\in \C$
be neither a critical value nor a locally omitted value.
Then there exists a simply connected region $D$ which does not
contain $a$, and such that $f^{-1}(D)$ is connected.
\end{theorem}

The case of a locally omitted value in Theorem~2
is based on the following result which is of independent interest.

\begin{theorem}\label{4} 
Let $f$ be an entire function of finite order,
and $a\in\C$ a locally omitted value. 
Then $a$ is the projection of a logarithmic
singularity of~$f^{-1}$.
\end{theorem}

The structure of the paper is the following.
In Section~2 we recall auxilliary facts on the singularities
of the inverses of entire functions.
In Section~3 we prove Theorems 1, 2 and 4.
In Section~4 we discuss some results needed for the 
proof of Theorem~3 and then we prove Theorem~3 in Section~5.
In Section~6 we analyse the example~(\ref{fexample}).  

\begin{ack}
We thank the referee for a very careful reading of the manuscript
and helpful suggestions, and in particular for pointing out an
error in the original argument.
\end{ack}

\section{Preliminaries}
\label{preliminaries}
We shall repeatedly use the following result of Iversen~\cite{I},
which
follows easily from the Gross Star Theorem~\cite[p.~292]{Ne},
or from the variant of the Gross Star Theorem stated
as Proposition~\ref{grossvariant} in Section 4 below.
\begin{ivth}
Let $\phi$ be a holomorphic branch of the inverse
$f^{-1}$ defined in a neighborhood of some point $w_0$ and let
$\gamma:[0,1]\to \C$ be a curve with $\gamma(0)=w_0$.
Then for every $\varepsilon>0$ there exists a curve $\tilde{\gamma}:[0,1]\to \C$
satisfying  $\tilde{\gamma}(0)=w_0$ and
$|\gamma(t)-\tilde{\gamma}(t)|<\varepsilon$  such that
$\phi$ has an analytic continuation  along~$\tilde{\gamma}$.
\end{ivth}

Now we recall Iversen's classification of singularities; 
see~\cite{AB}, \cite{I} or~\cite[p.~289]{Ne}.
Let $f$ be a transcendental meromorphic function and $a\in\C$.
Consider the open discs $B(a,r)$ of radius $r$ centered at~$a$.
For every $r>0$, it is possible to choose a component
$U_r$ of the preimage
$f^{-1}(B(a,r))$ in such a way that
$r_1<r_2$ implies $U_{r_1}\subset
U_{r_2}$. The possibility of such a choice
of (non-empty!) components $U_r$ follows
from Iversen's Theorem.

Now we have two possibilities:
\medskip

a) $\bigcap_{r>0}U_r$ consists of one point, or 
\medskip

b) $\bigcap_{r>0}U_r=\emptyset.$
\medskip

In the latter case we
say that our choice $r\mapsto U_r$ defines
a {\em transcendental singularity} of $f^{-1}$ over~$a$.
We also say that $a$ is the {\em projection}
of the transcendental singularity, or that the
transcendental singularity {\em lies over} $a$,
and any of the sets $U_r$ is called a {\em neighborhood}
of the transcendental singularity.
Projections of transcendental singularities coincide with
asymptotic values of~$f$.
A transcendental singularity over $a$ is called {\em direct} if for
some $r>0$ we have $f(z)\neq a$ for $z\in U_r$. Otherwise it is called
{\em indirect}. A direct singularity is called {\em logarithmic} if
the restriction $f:U_r\to B(a,r)\backslash\{ a\}$ is a universal covering
for some~$r>0$.
All these definitions can be also given for $a=\infty$ using 
$B(\infty,r)=\{ z\in\bC:|z|>1/r\}$.

It is clear from these definitions that locally omitted values
are exactly the projections of direct singularities.

For example, $\exp z$ has a logarithmic singularity over $0$,
and $(\sin z)/z$ has two indirect singularities over~$0$.

The importance  of direct singularities comes to a great extent from
the following result~\cite[\S XI.4]{Ne}.
\begin{dca}
A meromorphic function of finite order
has only finitely many direct singularities.
\end{dca}
A corollary of this result is that an entire function of finite order
has only finitely many asymptotic values~\cite[p.~313]{Ne}.

In Section~6, we will
prove that for the function (\ref{fexample})
the set of
direct singularities over $0$ has the power of continuum,
but none of these singularities 
is logarithmic. 
According to Heins~\cite{heins}, the set of projections of direct
singularities is always at most countable, but 
the set of direct singularities
over one point can have the power of the continuum.
Example~(\ref{fexample}) is a new example
of this kind; unlike the previous examples, it is given
by a simple explicit formula. 

\section{Proof of Theorems 1, 2 and 4}
{\em Proof of Theorem 1}. Suppose that $f^{-1}(D)$ is connected.
Using Iversen's Theorem, we can find a Jordan curve $\Gamma:[0,1]\to\C$
with $\Gamma(0)=\Gamma(1)=b$
for some point $b\in D$, such that $\Gamma$ does not pass through $a$,
$$\frac{1}{2\pi}\int_\Gamma\frac{dw}{w-a}=1,$$
and there exists a holomorphic branch
$\phi$ of $f^{-1}$ at $b$, such that $\phi$ has an analytic
continuation along~$\Gamma$. The preimage of $\Gamma$ under this branch
$\phi$ and its analytic continuation along $\Gamma$ is a simple compact
arc $\gamma$, which may be closed or not.
Both endpoints of $\gamma$ belong to $f^{-1}(D)$, and as $f^{-1}(D)$
is supposed to be connected,
we can find an arc $\gamma_1$ in $f^{-1}(D)$,
connecting the endpoints of~$\gamma$.
We have $f(\gamma')\subset D$, $D$ is simply connected, and $a\notin D$.
So 
$$\frac{1}{2\pi}\int_{f(\gamma')}\frac{dw}{w-a}=0.$$
So
$$\frac{1}{2\pi}\int_{\gamma\cup\gamma'}\frac{df(z)}{f(z)-a}=
\frac{1}{2\pi}\int_{\Gamma\cup f(\gamma')}\frac{dw}{w-a}=1,$$
which is a contradiction because $\gamma\cup\gamma'$ is a closed
curve and $f(z)\neq a$ in the plane.
This proves Theorem~1.
\medskip

Now we prove the following result from which Theorem~4 follows:

\begin{theorem}\label{theorem5} If an entire function has a direct
singularity over some point $a$ which
is not a logarithmic singularity, then every neighborhood 
of this singularity is also a neighborhood
of other direct singularities over~$a$.
\end{theorem}

It follows that whenever we have a direct
singularity over some point and no logarithmic singularities
over the same point, then the set of direct singularities
over this point has the power of the continuum.

As functions of finite order have only finitely
many direct singularities by the Denjoy-Carleman-Ahlfors
Theorem, we obtain Theorem~4.

\medskip

The proof of Theorem~5 requires the following lemma.

\begin{la} \label{poisson}
Let $\mu$ be a singular measure on the
unit circle, and $A=\{ e^{i\theta}:\theta\in (a,b)\}$ an
arc of the unit circle
such that $\mu(A)>0$. Then there exists
a point $\theta\in(a,b)$ such 
$$\lim_{r\to 1}u(re^{i\theta})=+\infty,$$
where  
$$u(re^{i\theta})=\frac{1}{2\pi}\int_{-\pi}^{\pi}
\frac{1-r^2}{1+r^2-2r\cos(t-\theta)}d\mu(t)$$
is the Poisson integral of~$\mu$.
\end{la}

{\em Proof.} For a subinterval $(x,y)$ of $(a,b)$, we denote by 
$\mu(x,y)$ the measure of the arc $\{ e^{i\theta}:\theta\in (x,y)\}$.
We first prove that there exists $\theta\in(a,b)$ such that
\begin{equation}
\label{phione}
\lim_{\varepsilon\to 0}\frac{1}{2\varepsilon}
\mu(\theta-\varepsilon,\theta+\varepsilon)=+\infty.
\end{equation}

Proving this by contradiction, suppose that such $\theta$ does not exist.
Then the sets
$$E_n=\{ x\in(a+1/n,b-1/n):\liminf_{\varepsilon\to 0}\frac{1}{2\varepsilon}
\mu(x-\varepsilon,x+\varepsilon)\leq n\},\quad n=n_0,n_0+1,\ldots$$
cover $(a,b)$.
Fix $n\geq n_0$. For every $x\in E_n$, 
there exists an interval of the form
$(x-\varepsilon,x+\varepsilon)$ with $\varepsilon\in(0,n^{-2})$
whose $\mu$-measure
is at most $4\varepsilon n$. By the well-known covering lemma
\cite[Thm. 1.1]{Guzman}, $E_n$ can be covered by some of these intervals
such that the multiplicity of this covering is an absolute
constant~$K$. Thus we obtain that $\mu(E_n)\leq 4K/n.$
As the sets $E_n$ form an increasing sequence,
we conclude that $\mu(E_n)=0$ for all $n\geq n_0$.
So $\mu(a,b)=0$ and we obtain a contradiction,
which proves the existence of the point $\theta$ satisfying~(\ref{phione}).

Now it is easy to pass from (\ref{phione}) to the Poisson integral.
For $0<\varepsilon<\pi$ we have
\begin{eqnarray*}
u(re^{i\theta})
&\geq&
\frac{1}{2\pi}\int_{\theta-\varepsilon}^{\theta+\varepsilon}
\frac{1-r^2}{1+r^2-2r\cos(t-\theta)}d\mu(t)\\
&\geq&
\frac{1}{2\pi} \frac{1-r^2}{1+r^2-2r\cos\varepsilon} 
\mu(\theta-\varepsilon,\theta+\varepsilon).
\end{eqnarray*}
Putting $r=1-\varepsilon$ and noting that
then 
$$\frac{1-r^2}{1+r^2-2r\cos\varepsilon}\geq \frac{1}{\varepsilon}$$
for sufficiently small $\varepsilon$ 
we obtain
$$u(re^{i\theta})\geq
\frac{1}{2\pi\varepsilon} 
\mu(\theta-\varepsilon,\theta+\varepsilon),$$
which completes the proof.

\medskip

{\em Proof of Theorem 5.}
Suppose that $U=U_r$ is a neighborhood of
exactly one direct
singularity over a finite point $a$, where $r>0$
is so small that $f(z)\neq a$ in~$U$. We are going to prove that
this singularity is logarithmic.

By the Maximum Principle, $U$ is simply connected.
It is easy to see
that the closure of $U$ in the Riemann sphere
is locally connected.
So a conformal map $\phi:B(0,1)\to U$ extends to a continuous
map from the unit disc to the Riemann sphere.
The preimage of infinity under $\phi$ is a closed subset $E$,
of the unit circle,
which by a theorem of
Beurling~\cite[p.~344]{Tsuji} has zero logarithmic capacity.

We consider the positive harmonic function
$$u(z):=\log\frac{r}{|f(\phi(z))-a|},\quad z\in B(0,1).$$
It has a Poisson representation
$$u(re^{i\theta})=
\frac{1}{2\pi}\int_{-\pi}^{\pi}
\frac{1-r^2}{1+r^2-2r \cos(t-\theta)} d\mu(t)$$
for some finite Borel measure $\mu$, and we
have $\lim_{r\to 1} u(re^{i\theta})=0$ if
$e^{i\theta}\notin E$. As $E$ has zero capacity,
and thus zero length, the measure $\mu$ is singular.
It is easy to see that 
if $\mu$ is a single atom so that $u$ is proportional
to the Poisson kernel, then 
$f:U\to B(r,a)\setminus \{a\}$ is a universal covering, 
and thus the singularity
we consider is logarithmic.

Otherwise, there is a simple cross-cut $\sigma$ in $B(0,1)$
beginning and ending in the complement of $E$,
such that the two arcs on the unit circle bounded by 
the endpoints of  $\sigma$ both intersect the support of~$\mu$.
Lemma~\ref{poisson} implies that 
$u$ is unbounded in each of the two components
$G_1$ and $G_2$ of $B(0,1)\backslash\sigma$.
The image $\phi(\sigma)$
of this cross-cut separates $U$ into two regions $D_j=\
\phi(G_j)$.
The harmonic function
$$v(z):=u(\phi^{-1}(z))=\log\frac{r}{|f(z)-a|}$$
is unbounded in each $D_j$ and bounded on $\partial D_j$ for $j=1,2$.
Thus there exists $\varepsilon>0$ such that
$$\left\{ z\in U: v(z)>\log\frac{r}{\varepsilon}\right\}=
\{ z\in U:|f(z)-a|<\varepsilon\}$$
is disconnected. As $f(z)\neq a$ for $z\in U$ we conclude that
$U$ is a neighborhood of at least two singularities over~$0$.
This completes the proof of Theorem~5.
\medskip

As we already mentioned, Theorem~4 follows from Theorem~5 and
the Denjoy-Carleman-Ahlfors Theorem.

Now Theorem~2 is an easy corollary:
In the case of a critical value, we repeat Baker's argument~\cite{Ba}
and in the case of a locally omitted value, we first use Theorem~4,
to conclude that this singularity is in fact logarithmic,
and then repeat the argument from~\cite{EL}.
Both~\cite{Ba} and~\cite{EL} deal only with the case that $D$ is 
completely invariant, but the arguments extend to the 
situation of Theorem~2 without difficulty.

\section{Results needed for the proof of Theorem~3}
The following definition will be used in the proof of
Theorem~3. 
Let $f$ be an entire function.
A simple curve 
$\gamma$ 
will be called {\em good for $f$},
if $\gamma$ contains
no critical values of $f$ and all
components of the full preimage $f^{-1}(\gamma)$
are compact.

It is easy to see that a simple curve
which contains neither critical values nor asymptotic values
is good. For entire functions of finite order,
there can be only finitely many asymptotic values by 
the Denjoy-Carleman-Ahlfors Theorem. Thus
we obtain the existence of good curves, and in
fact we see that the conclusion of Proposition~\ref{grossvariant} 
and~\ref{goodinterval} below holds for entire functions of finite order.

In general, there are entire functions
for which every point in the complex plane is an
asymptotic value~\cite{Gross}, so the existence of good curves for such
functions is not evident.  
An instructive example is given by 
$f(z) = (\sin z)/z$. Here $0$ is the projection of an indirect
singularity. However, one can show that the segment $[-i,i]$ is
good. On the other hand, $[-\varepsilon,\varepsilon]$ is not good for
any positive $\varepsilon$.

In the remaining part of this section we will
prove the existence of good curves in general.
This material is not used anywhere else in the paper, but may be
of independent interest.

Our Proposition~\ref{grossvariant} and~\ref{goodinterval}
below are similar to the results of
Shimizu~\cite[p.~186]{Shi} and
Terasaka~\cite[Lemma on p.~310]{Ter}.
We need the following version
of the classical Gross Star Theorem~\cite[p.~292]{Ne}:

\begin{prop}
\label{grossvariant}
Let $\phi$ be a holomorphic branch of the inverse $f^{-1}$ defined
in some disc~$B$, and let $\ell$ be some direction in the plane.
Then $\phi$ has an analytic continuation along almost all straight
lines intersecting $B$ and having the direction~$\ell$.
\end{prop}

Such lines can be parametrized by the points of their intersection
with the diameter of $B$ perpendicular to~$\ell$. ``Almost all''
refers to the Lebesgue measure on this diameter.
\medskip

{\em Proof of Proposition~\ref{grossvariant}.} We assume for simplicity that
the direction $\ell$ is parallel to the real axis, and
that the diameter of $B$ perpendicular to this direction is $(ia,ib)$,
where $a<b$.
Let $M>b-a$.
Consider the rectangle
$$Q_M:=\{ z=x+iy: |x|<M, y\in (a,b)\}.$$
For each horizontal interval $\{ x+iy_0: |x|<M\}$, where $y_0\in(a,b)$,
we consider the maximal open subinterval containing the point $iy_0$
such that an analytic continuation of $\phi$ is possible
along this subinterval.
The union of these maximal subintervals over all $y_0\in(a,b)$
forms a region $G_M\subset Q_M$.
If a maximal horizontal interval in $G_M$
has an endpoint
inside $Q_M$,
then we will call this endpoint a {\em singular
point of}~$\phi$. It is enough to
show that the Lebesgue measure
of the projection of the set of singular points
on the imaginary axis is zero,
for every fixed $M>b-a$.
The analytic continuation of $\phi$
along maximal horizontal intervals
in $G_M$ maps
$G_M$ univalently onto some region $G_M^\prime\subset\C$.
The singular points
in $G_M$ correspond
to the critical values of $f$ and
to the accessible points at infinity of~$G_M^\prime$.
Since the set of critical values is countable,
the Lebesgue measure
of its projection
on the imaginary axis is zero.

Let $\sigma^\prime_r$ be
the intersection of $G_M^\prime$ with the circle
$\{ z:|z|=r\}$.
We may assume that $\phi$ is bounded in~$B$.
Then for $r>r_0$
the set $\sigma_r:=f(\sigma^\prime_r)$
is a union of cross-cuts
in $G_M$ which separate the diameter $[ia,ib]$
from the set of singular points
of $\phi$ on the boundary of~$G_M$.
It is enough to show that the
length of $\sigma_r$ tends to zero
as $r\to\infty$ on some sequence.

We have
$$\length(\sigma_r)=\int_{\sigma_r^\prime}|f'(z)|\,
|dz|$$
and by Schwarz's inequality
$$\length^2(\sigma_r)\leq 2\pi r
\int_{\sigma_r^\prime}|f'(z)|^2\, |dz|.$$
Dividing by $r$ and integrating with respect to $r$
from $r_0$
to $\infty$, we obtain
$$\int_{r_0}^\infty\length^2(\sigma_r)\frac{dr}{r}\leq
2\pi\int\!\!\int_{G_M'}|f'(z)|^2\, dxdy=2\pi \area(G_M)\leq
2\pi  \area(Q_M).$$
Thus the integral on the left hand side converges.

We conclude that $\length(\sigma_r)\to 0$
on some sequence $r=r_k\to\infty$. This proves the Proposition.

\begin{prop}
\label{goodinterval}
Let $f$ be an entire function, and
$Q=(a,b,c,d)$ a rectangle in the plane.
Then almost every closed interval
connecting the opposite
sides $[a,b]$ and $[c,d]$ and parallel
to the other sides is good for~$f$.
\end{prop}

{\em Proof}. Consider the set of pairs $\{ B_j,\phi_j\}$,
where $B_j$ is a disc contained in $Q$,
having rational center and rational radius, and
$\phi$ is a holomorphic branch of $f^{-1}$
in this disc. According to the Poincar\'e--Volterra
theorem, this set is countable.
Applying
Proposition~\ref{grossvariant} to the pair $\{ B_j,\phi_j\}$
and the
direction $[a,d]$ we obtain an exceptional set of lines $E_j$
of measure zero.
Then $E=\bigcup E_j$ is a set of measure zero,
and all intervals which are intersections of $Q$ with
lines parallel to $[a,d]$ and not in $E$ are good for~$f$. This proves the
Proposition.
\medskip

As mentioned, we will use Propositions~\ref{grossvariant} and~\ref{goodinterval}
only for entire functions of finite order, and for such functions the
conclusion follows from the Denjoy-Carleman-Ahlfors Theorem.
\section{Proof of Theorem 3}
\begin{la}
\label{connected}
Let $f$ be an entire function of finite order and let $a\in\C$.
Then there exists $r_0>0$ such that if $0<r\leq r_0$, then
all components of
$f^{-1}(B(a,r))$
have connected boundary.
\end{la}
{\em Proof}.
Let $U$ be a component of $f^{-1}(B(a,r))$.
By the Maximum Principle, $U$ is simply connected.
Each complementary component of $U$
contains a neighborhood of a singularity of $f^{-1}$ over~$\infty$.

By the Denjoy-Carleman-Ahlfors Theorem, there are only finitely
many singularities of  $f^{-1}$ over~$\infty$.
It follows that there exists $r_0>0$ such that
$f^{-1}(\C\setminus\overline{B(a,r)})$ is connected for $0<r\leq r_0$.
As this set is a neighborhood of each singularity of $f^{-1}$ over $\infty$
we deduce that the complement of $U$ is connected. Hence the boundary
of $U$ is connected.
\medskip

{\em Proof of Theorem 3}.
We will construct a simple curve $\Gamma$ connecting $a$ with $\infty$
such that the preimage of $\Gamma$ consists of infinitely many simple,
pairwise disjoint
curves which connect the preimages of $a$ with~$\infty$.
With $D:=\C\backslash\Gamma$ we then see that both
$D$ and $f^{-1}(D)$ are simply connected.

We choose $r_0$ according to Lemma~\ref{connected}.
We begin with a simple curve $\Gamma^0$ (for example,
a straight line segment)
which is good for $f$ and connects a point $w_0\in B(a,r_0)$ to
some point in $\C\backslash\overline{B(a,r_0)}$.
In addition we may
assume that $\Gamma^0\cap \partial B(a,r)$ consists
of at most one point for all~$r$. Curves with the latter property
will be called {\em $a$-monotonic}.
The existence of a curve $\Gamma^0$ with the properties mentioned
follows from Proposition~\ref{goodinterval}.

Take a point
a point $b\in \Gamma^0\cap \C\backslash\overline{B(a,r_0)}$
which is not the endpoint of~$\Gamma^0$ and
let $(c_j)$ be the sequence of $b$-points of~$f$.
Let $\gamma_1^0$ be the component of
$f^{-1}(\Gamma^0)$ which contains~$c_1$.
Then $\gamma_1^0$ is a simple curve connecting a $w_0$-point
$x_1$ with~$c_1$.
Let $U_1$ be the component of $f^{-1}(B(a,r_0))$ that contains~$x_1$.
Since $f$ has no direct singularity or critical point
over $a$, there exists $z_1\in U_1$ with $f(z_1)=a$ and
$f'(z_1)\neq 0$.
Thus there exists $r_1$ with $0<r_1<r_0$ such that
there is a branch $\phi_1$ of $f^{-1}$ which is
defined in $B(a,r_1)$ and maps $a$ to~$z_1$.
We may also assume that $\phi_1$  is bounded in  $B(a,r_1)$.

We can connect $z_1$ by a curve $\sigma_1$ to $\partial
U_1$  such that $f(\sigma_1)$ is a straight line connecting
$a$ to $\partial B(a,r_0)$.
(Here we say that a curve $\gamma$ {\em connects} a point
$z$ to a set $S$ if $\gamma$ is a simple curve such that one
endpoint of $\gamma$ is $z$ while the other one is in $S$,
and $S\cap \gamma$ consists only of that endpoint.)
Since $r_0$  has been chosen according to Lemma~\ref{connected},
the boundary of $U_1$ is connected.
Thus we can connect the endpoint of $\sigma_1$ in $\partial U_1$
by a curve $\sigma'_1\subset \partial U_1$
to the point $v_1$ which lies in the intersection
of $\gamma_1^0$ and~$\partial U_1$.
(Note that the intersection of $\gamma_1^0$ and $\partial U_1$
consists of only one point since $\Gamma^0=f(\gamma_1^0)$ is
$a$-monotonic.)
The curve  $\sigma_1+\sigma'_1$  thus connects $z_1$ to~$v_1$.
By deforming $\sigma'_1$ slightly we can replace
the curve  $\sigma_1+\sigma'_1$  by a curve
$\tau_1$ which connects $z_1$ to $v_1$ such that
$f(\tau_1)$ is $a$-monotonic.
Using Proposition~\ref{goodinterval}, we can replace $\tau_1$ by a
curve $\tau'_1$  which connects a point $y_1^1\in \phi_1(B(a,r_1))$
to $\gamma_1^0\cap U_1$ and which has the property that $f(\tau'_1)$ is good
and $a$-monotonic.
Combining $f(\tau'_1)$ and $\Gamma^0$ we thus obtain a curve
$\Gamma^1$ which is good
and $a$-monotonic and which connects $w_1:=f(y_1^1)\in B(a,r_1)$
to~$b$.
More precisely, if $u_1$ is the endpoint of $f(\tau'_1)$ in
$\Gamma^0$ and if $\Sigma^0$ is the arc that connects $u_1$
and $b$ in $\Gamma^0$, then we take $\Gamma^1:=f(\tau'_1)\cup \Sigma^0$.
Note that $u_1\in B(a,r_0)$ since the endpoint of $\tau'_1$
is in~$U_1$.

The component $\gamma_1^1$ of $f^{-1}(\Gamma^1)$ that
contains $c_1$ consists of $\tau'_1$ and a subarc of $\gamma_1^0$,
and it is a simple curve
connecting $y_1^1\in \phi(B(a,r_1))$ to~$c_1$.
In fact, since $\Gamma^1$ is good, we
see that for all for $j\in \N$ the component
$\gamma_j^1$ of $f^{-1}(\Gamma^1)$ which contains $c_j$
is a simple curve connecting a $w_1$-point $y_j^1$ to~$c_j$.

The following fact is important: {\em
no matter how we
extend $\Gamma^1$ by attaching a piece
in $B(a,r_1)$,
the component of the preimage of the extended curve
that contains $c_1$ and hence $\gamma^1_1$ will be compact}.
This follows since
the part added to $\gamma^1_1$ will be
contained in $\phi_1(B(a,r_1))$.

Now we repeat this process of extension.
Suppose that
$r_{n-1}<r_{n-2}<\dots<r_1<r_0$ and that
$\Gamma^{n-1}$ is a good and $a$-monotonic curve
which connects a point $w_{n-1}\in B(a,r_{n-1})$ to $b$ such
that for $1\leq j\leq {n-1}$ the component $\gamma_j^{n-1}$ of
$f^{-1}(\Gamma^{n-1})$ which contains $c_j$ has the following
property:
{\em no matter how we extend $\Gamma^{n-1}$ by attaching a piece
in $B(a,r_{n-1})$,
the component of the preimage of the extended curve
that contains $c_j$ and hence $\gamma^1_j$ will be compact
for $1\leq j\leq {n-1}$}.

The way we obtain $\Gamma^n$ from $\Gamma^{n-1}$ is essentially the same
that we used to obtain $\Gamma^1$ from $\Gamma^0$:
since $\Gamma^{n-1}$ is good, the component
$\gamma_n^{n-1}$ of $f^{-1}(\Gamma^{n-1})$
that contains $c_{n}$ is a simple curve connecting a $w_n$-point $x_n$
to~$c_n$.
Let $U_{n}$ be the component of $f^{-1}(B(a,r_{n-1}))$ that contains~$x_n$.
Then
there exists $z_n\in U_n$ with $f(z_n)=a$ and
$f'(z_n)\neq 0$, and hence
there exists $r_n$ with $0<r_n<r_{n-1}$ such that
there is a branch $\phi_n$ of $f^{-1}$ in $B(a,r_n)$
with $\phi_n(a)=z_n$, where we
may again assume that $\phi_n$  is bounded in  $B(a,r_n)$.
We connect $z_n$ by a curve $\sigma_n$ to $\partial U_n$
such that $f(\sigma_n)$ is a straight line
and we connect the  endpoint of $\sigma_n$
by a curve $\sigma'_n\subset \partial U_n$
to the point $v_n$ which lies in the intersection
of $\gamma_n^{n-1}$ and~$\partial U_n$.
Again we can replace
the curve  $\sigma_n+\sigma'_n$  by a curve
$\tau_n$ which connects $z_n$ to $v_n$ such that
$f(\tau_n)$ is $a$-monotonic and,
using Proposition~\ref{goodinterval}, we can replace $\tau_n$ by a
curve $\tau'_n$  which connects a point $y_n^n\in \phi_n(B(a,r_n))$
to $\gamma_n^{n-1}\cap U_n$
and which has the property that $f(\tau'_n)$ is good
and $a$-monotonic. From
$f(\tau'_n)$ and $\Gamma^{n-1}$ we now obtain a curve
$\Gamma^n$ which is good
and $a$-monotonic and which connects $w_{n}:=f(y_n^n)\in B(a,r_n)$
to~$b$.
Moreover, the component $\gamma_n^n$ of $f^{-1}(\Gamma^n)$ that
contains $c_n$ is a simple curve
connecting $y_n^n\in \phi_n(B(a,r_n))$ to $c_n$,
and no matter how we extend
$\Gamma^n$ by attaching a piece in $B(a,r_n)$,
the component $\gamma_n^n$ of the preimage of the extended curve
which contains $c_n$ will be compact.
And it follows from our induction hypothesis that
the same is true for the preimages $\gamma_j^n$ of the extended curve
which contain $c_j$, for $1\leq j\leq n-1$.

Note that $\Gamma_n$ need not contain $\Gamma_{n-1}$,
but since the endpoint of $f(\tau'_n)$ is contained in
$B(a,r_{n-1})$ we have
$$\Gamma_n\setminus B(a,r_{n-1})= \Gamma_{n-1}\setminus
B(a,r_{n-1})\supset \Gamma_{n-1}\setminus B(a,r_{n-2})$$
for $n\geq 2$.

We now combine the curves $\Gamma_n$ defined inductively in the above
way to a curve $\Gamma^\infty$ by putting
$$\Gamma^\infty:=\{a\}\cup
\bigcup_{n=1}^\infty \left(\Gamma_n\setminus B(a,r_{n-1})\right).$$
Then $\Gamma^\infty$ is a simple (and in fact $a$-monotonic)
curve that connects $a$ to $b$, and
it follows from the construction of the $\Gamma_n$ that
the preimage of $\Gamma^\infty$ that contains $c_j$ is a simple
curve that connects $z_j$ with~$c_j$.

Finally we  connect $b$ to $\infty$ by an $a$-monotonic curve
$\Sigma$ with the property that every compact subarc of
$\Sigma$  is good. Such a curve exists by Proposition~\ref{goodinterval}.
It then follows that
$\Gamma:=\Gamma^\infty\cup \Sigma$ has the properties stated
at the beginning of the proof.

This completes the proof of Theorem~3.

\section{An example}\label{example}
We show that the function $f$ given 
by (\ref{fexample})
has infinitely many direct but no  logarithmic singularity over~$0$.
Let 
$$g(z):=\sum_{k=1}^\infty 
\left(\frac{z}{2^k}\right)^{2^k}$$
so that $f(z)=\exp g(z)$.
We fix $\varepsilon$ with $0<\varepsilon\leq\frac18$ and 
put $r_n:=(1+\varepsilon) 2^{n+1}$ and 
$r_n':=(1-2\varepsilon)2^{n+2}$ for 
$n\in\N$. For $j\in\{0,1,\dots,2^n-1\}$ we define the sets
$$A_{j,n}:=
\left\{r\exp\left(\frac{2\pi i j}{2^n}\right): r\geq 
r_n\right\},\quad i=\sqrt{-1},$$
$$B_{j,n}:=
\left\{r\exp\left(\frac{\pi i}{2^n}+
\frac{2\pi i j}{2^n}\right): 
r_n\leq r\leq r_n'\right\},$$ 
and 
$$C_{j,n}^\pm:= 
\left\{r\exp\left(\frac{\pi i}{2^n}+ \frac{2\pi i j}{2^n} \pm 
\frac{r- r_n'}{r_{n+1}-r_n'} 
\frac{\pi i}{2^{n+1}}\right): 
r_n'\leq r\leq r_{n+1}\right\}.$$ 
We shall show that if $n$ is large enough, then 
\begin{equation}\label{reallarge}
\re g(z) > 2^{2^n} \quad \text{for}\quad z\in A_{j,n}
\end{equation}
while
\begin{equation}\label{realsmall}
\re g(z) < -2^{2^n} \quad \text{for}\quad z\in B_{j,n}\cup C^+_{j,n}
\cup C^-_{j,n}.
\end{equation}
Note that $C^-_{j,n}$ connects $B_{j,n}$ to $B_{2j,n+1}$ 
while $C^+_{j,n}$ connects $B_{j,n}$ to $B_{2j+1,n+1}$.
This implies that 
$$T:= \left[-ir_1,ir_1\right]\cup \bigcup_{n=1}^\infty
\bigcup_{j=0}^{2^n-1}
\left(B_{j,n} \cup C^+_{j,n}  \cup C^-_{j,n}
\right)
$$
is an
infinite binary tree; see Figure~1.
\begin{figure}[htb]
\begin{center}
\includegraphics[width=26pc]{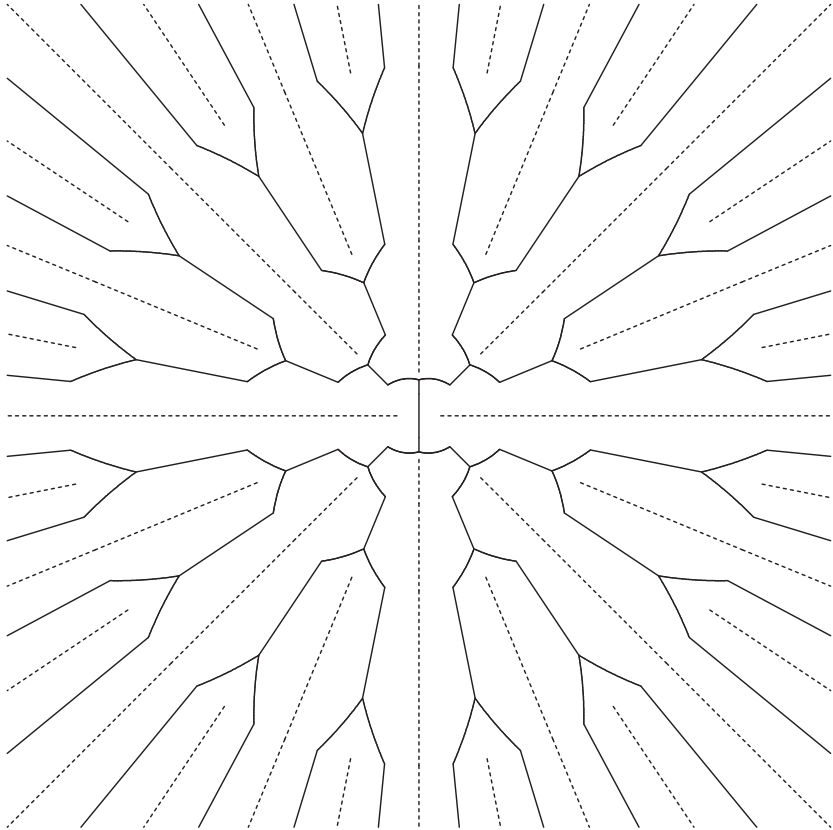}
\caption{The part of the tree $T$ lying in $\{z:|\re z|\leq 80, 
|\im z|\leq 80\}$, for $\varepsilon=1/16$.
The sets $A_{j,n}$ are drawn as dotted lines.}
\end{center}
\end{figure}
By (\ref{realsmall}), every unbounded simple path 
on this tree starting at $0$ 
is an asymptotic curve on which $\re g(z)\to -\infty$.
Choosing $U_\rho$ as the component of 
$$\{z: \re g(z)< \log \rho\}= \{z: |f(z)|< \rho\}$$
which contains 
the ``tail'' of this curve we thus obtain a 
transcendental singularity  of $f^{-1}$ over $0$,
and this singularity is direct because $f$ has no zeros.
Using (\ref{reallarge}) we see that 
different curves define different singularities.
Thus we obtain a set of direct singularities which has the
power of the continuum. 

Moreover, it follows from (\ref{reallarge}) and (\ref{realsmall})
and the above considerations 
that if $U_\rho$ is a component of
$\{z: |f(z)|< \rho\}$ 
containing the ``tail'' of some curve in $T$, then
$U_\rho$ also contains the ``tail'' of some  other curve 
in $T$ and thus 
there exists $\rho'<\rho$ such that $U_\rho$ contains 
at least two components of $\{z: |f(z)|< \rho'\}$.
This implies that the singularity defined by $\rho\mapsto U_\rho$
is not logarithmic.

To prove (\ref{reallarge}) we note that if 
$z=r\exp\left(2\pi i j/2^n \right)\in A_{j,n}$ 
so that 
$r\geq r_n$,
then
$$\re g(z)\geq 
\sum_{k=n}^\infty 
\left(\frac{r}{2^k}\right)^{2^k}
- \sum_{k=1}^{n-1}
\left(\frac{r}{2^k}\right)^{2^k}
\geq 
\left(\frac{r}{2^n}\right)^{2^n}
- \sum_{k=1}^{n-1}
\left(\frac{r}{2^k}\right)^{2^k}.
$$
Put $s:=r/2^n$ and $$\Sigma_1:=\sum_{k=1}^{n-1}
\left(\frac{r}{2^k}\right)^{2^k}.$$
Then 
$$\Sigma_1=\sum_{k=1}^{n-1} \left(s2^{n-k}\right)^{2^k}
\leq s^{2^{n-1}}\sum_{k=1}^{n-1}2^{(n-k)2^k}.$$
Now $(n-k)2^k\leq 2^{n-1}$ for $1\leq k\leq n-1$ and $s<2+2\varepsilon$
so that
$$\Sigma_1\leq s^{2^{n-1}}(n-1)2^{2^{n-1}}
=o\left( s^{2^{n}}\right)$$
as $n\to\infty$ 
and hence
$$\re g(z)\geq
\left(\frac{r}{2^n}\right)^{2^n}-\Sigma_1
= (1-o(1))s^{2^{n}}> 2^{2^{n}}$$
for large~$n$.
To prove (\ref{realsmall}) for $z\in B_{j,n}$,
let 
$z=r\exp\left(\pi i/2^n+2\pi i j/2^n \right)\in B_{j,n}$,
with
$r_n\leq r\leq r_n'$.
Then 
$$\re g(z)\leq 
-\left(\frac{r}{2^n}\right)^{2^n}+\Sigma_1
+\Sigma_2$$
with 
$$\Sigma_2
:=\sum_{k=n+1}^\infty 
\left(\frac{r}{2^k}\right)^{2^k}
=\sum_{k=n+1}^\infty  \left(s2^{n-k}\right)^{2^k}.$$
Thus 
$$\Sigma_2\leq 
\left(\frac{s}{2}\right)^{2^n}
+\sum_{k=n+2}^\infty
\left(\frac{s}{4}\right)^{2^k}.$$
Since $s/4\leq 1-2\varepsilon$ we find that 
$$\Sigma_2
=o\left(s^{2^{n}}\right)$$
as $n\to\infty$
and thus 
$$\re g(z)\leq -(1-o(1))s^{2^{n}}<-2^{2^{n}}$$
for $z\in B_{j,n}$, provided $n$ is sufficiently large.

Finally we prove (\ref{realsmall}) for $z\in C^+_{j,n}$.
So let 
\begin{eqnarray*}
z & = &    
r\exp\left(\frac{\pi i}{2^n}+ \frac{2\pi i j}{2^n}+
\frac{r- r_n'}{r_{n+1}-r_n'} 
\frac{\pi i}{2^{n+1}}\right)
\\
& = & 
r\exp\left(\frac{\pi i}{2^n}+ \frac{2\pi i j}{2^n}+
\frac{s-4(1-2\varepsilon)}{12\varepsilon}
\frac{\pi i}{2^{n+1}}\right)\\
& \in &
C^+_{j,n},
\end{eqnarray*}
with 
$r_n'\leq r\leq r_{n+1}$,
so that $4(1-2\varepsilon)\leq s \leq 4(1+\varepsilon).$
We have 
$$\re g(z)\leq 
\re\left(\frac{z}{2^n}\right)^{2^n}+
\re\left(\frac{z}{2^{n+1}}\right)^{2^{n+1}}+
\Sigma_1
+\Sigma_3$$
with 
$$\Sigma_3
:=\sum_{k=n+2}^\infty 
\left(\frac{r}{2^k}\right)^{2^k}
\leq \sum_{k=n+2}^\infty 
\left(\frac{s}{4}\right)^{2^k}=o(1)$$
since $s>4$.
Since $\Sigma_1=o\left(s^{2^n}\right)$ 
we find  that 
\begin{eqnarray*}
\re g(z)
&\leq&
 s^{2^n}
 \cos\left(\pi +\frac{s-4(1-2\varepsilon)}{12\varepsilon}
 \frac{\pi}{2}\right)\\
& & 
 +\left(\frac{s}{2}\right)^{2^{n+1}}
\cos\left(\frac{s-4(1-2\varepsilon)}{12\varepsilon}\pi\right)
+o\left( s^{2^n}\right)\\
&=&
 s^{2^n}
 \left(
 \cos\left(\pi+t \frac{\pi}{2}\right)
 +\left(\frac{s}{4}\right)^{2^{n}}
 \cos\left(t\pi\right) +o(1) \right)
\end{eqnarray*}
as $n\to \infty$, 
with $t:=(s-4(1-2\varepsilon))/12\varepsilon$.
The range $4(1-2\varepsilon)\leq s\leq 4(1+\varepsilon)$
corresponds to $0\leq t\leq 1$ and 
$s=4(1-2\varepsilon)+12\varepsilon t$.
We define
\begin{eqnarray*}
h(t)
&:=&
\cos\left(\pi+t \frac{\pi}{2}\right)
 +\left(\frac{s}{4}\right)^{2^{n}}
  \cos\left(t\pi\right)\\
&=&
\cos\left(\left(1+\frac{t}{2}\right)\pi\right)
+(1-2\varepsilon+3\varepsilon t)^{2^n}\cos\left(t\pi\right).
\end{eqnarray*}
and put $\delta:=-\cos(11\pi/8)/2>0$.
For $0\leq t\leq \frac12$ we have 
$$h(t)\leq \cos\left(\frac54 \pi\right)+\left(1-\frac{\varepsilon}{2}\right)^{2^{n}}
<-2\delta$$
if $n$ is large enough.
For $\frac12\leq t\leq  \frac34$ we have 
$\cos\left(t\pi\right)<0$ and thus
$$h(t)\leq \cos\left(\left(1+\frac{t}{2}\right)\pi\right)
\leq  \cos\left(\frac{11}{8}\pi\right)=-2\delta.$$
Finally, for $\frac34\leq t\leq  1$ we have 
$\cos\left(\left(1+\frac{t}{2}\right)\pi\right)<0$ and thus
$$h(t)\leq (1-2\varepsilon+3\varepsilon t)^{2^n}\cos\left(t\pi\right)
\leq \left(1+\frac{1}{4}\right)^{2^{n}}\cos\left(\frac{3}{4}\pi\right)
\leq -2\delta$$
if $n$ is large.
Overall we find that 
$h(t)\leq  -2\delta$ for all $t$ and thus
$$ \re g(z) \leq -\delta s^{2^n} < - 2^{2^n}$$
for $z\in C^+_{j,n}$, provided $n$ is large enough.
The proof that 
$$ \re g(z) < - 2^{2^n}$$
for $z\in C^-_{j,n}$ is analogous.
This completes the proof of (\ref{reallarge}) and
(\ref{realsmall}). As already mentioned, this implies 
that  every path going
to $\infty$ in $T$ corresponds 
to a direct singularity of $f$ over $0$ 
which is not logarithmic, 
and the set of such singularities has the power of 
the continuum.
Also,
we see that if $\rho\to U_\rho$  is a singularity 
over $0$
such that $U_\rho\cap T\neq \emptyset$ for all $\rho>0$,
then this singularity is not logarithmic.

It remains to prove that there are no other 
singularities over~$0$.
Suppose that $\rho\to U_\rho$ 
is a singularity
over $0$
such that $U_\rho\cap T= \emptyset$ for some $\rho>0$.
In order to obtain a contradiction 
we note that 
it follows as in the proof of (\ref{reallarge}) and (\ref{realsmall})
that if 
$r_n\leq |z|\leq r_n'$,
then 
$$g(z)=(1+\eta(z))
\left(\frac{z}{2^n}\right)^{2^n}$$
where $\eta(z)\to 0$ as $n\to \infty$.
For large $n$ we thus have $|\eta(z)|\leq \varepsilon^2\leq \frac12$.
Differentiating we obtain
$$\frac{g'(z)}{g(z)}-\frac{2^n}{z}=\frac{\eta'(z)}{1+\eta(z)}.$$
For 
$(1+2\varepsilon) 2^{n+1}\leq |z|\leq (1-3\varepsilon)2^{n+2}$
we thus have 
\begin{eqnarray*}
\left|\frac{g'(z)}{g(z)}-\frac{2^n}{z}\right|
&\leq& 
2 |\eta'(z)|\\
&=& 
2 \left|\frac{1}{2\pi i}\int_{|\zeta-z|=\varepsilon 2^{n+1}}
\frac{\eta(\zeta)}{(\zeta-z)^2}d\zeta \right|\\
&\leq&
2 \frac{1}{\varepsilon 2^{n+1}} \max_{|\zeta-z|=\varepsilon 2^{n+1}}
|\eta(z)|\\
&\leq&
\frac{\varepsilon}{2^{n}} 
\end{eqnarray*}
and hence 
$$\left|\frac{z g'(z)}{g(z)}-{2^n}\right| \leq \frac{\varepsilon |z|}{2^{n}} 
 \leq 4 \varepsilon (1-3\varepsilon) <\frac12.$$
We deduce that 
\begin{eqnarray*}
\frac{ d\arg g(re^{i\theta})}{d\theta}
&=&
\im \left(\frac{ d\log g(re^{i\theta})}{d\theta}\right)\\
&=&
\re \left( \frac{re^{i\theta}g'(re^{i\theta})}{g(re^{i\theta})} \right)\\
&\geq &
2^{n}- \frac12 \\
&> & 0
\end{eqnarray*}
for $(1+2\varepsilon) 2^{n+1}\leq r\leq (1-3\varepsilon)2^{n+2}$
and large~$n$. 
We conclude that  $\arg g(re^{i\theta})$ is an increasing function
of $\theta$, and it increases by $2^n 2\pi$
as $\theta$ increases by~$2\pi$. 
Choose $n$ and $r$ as above so large that 
the circle $\{z:|z|=r\}$
intersects $U_\rho$, that~(\ref{reallarge}) and~(\ref{realsmall}) 
hold and that $-2^{2^n}<\log \rho$. 
From the behavior of 
$\arg g(re^{i\theta})$ we deduce that 
the circle
$\{z:|z|=r\}$ contains at most $2^n$ arcs where 
$\re g(re^{i\theta})<\log \rho$.
On the other hand, for $j\in\{0,1,\dots,2^n-1\}$ the points 
$r\exp\left(\pi i/2^n+2\pi i j/2^n\right)$
are contained in such an arc by~(\ref{realsmall}),
and each of them is  contained in  a different one
by~(\ref{reallarge}).
Hence there are precisely $2^n$ such arcs and each one
contains one of the points 
$r\exp\left(\pi i/2^n+2\pi i j/2^n\right)$. Thus 
each such arc intersects some $B_{j,n}$ and hence~$T$.
In particular, $U_\rho\cap \{z:|z|=r\}$ intersects 
$T$, contradicting the assumption that $U_\rho\cap T
=\emptyset$. 
This completes the proof that 
$f$ has no  logarithmic singularities over~$0$.
\medskip

\begin{rem}
It is much easier to find meromorphic functions with 
a direct singularity which is not logarithmic.
For example, $f(z)=1/(z \sin z)$ has two direct singularities 
over $0$, but none of them is logarithmic, since their neighborhoods
are multiply connected. 
\end{rem}

\end{document}